\documentclass{amsart}

\usepackage{amssymb}
\usepackage{amsmath}
\usepackage{amsthm}
\usepackage{amsbsy}
\usepackage[dvips]{graphics}
\usepackage{amscd}
\usepackage{bm}

\theoremstyle{plain}
\newtheorem{theorem}{Theorem}[section]

\newtheorem{lemma}[theorem]{Lemma}
\newtheorem{corollary}[theorem]{Corollary}
\newtheorem{proposition}[theorem]{Proposition}

\theoremstyle{definition}
\newtheorem{definition}{Definition}[section]

\theoremstyle{remark}
\newtheorem{remark}[theorem]{Remark}
\newtheorem*{acknowledgements}{Acknowledgements}
\newtheorem{example}{Example}[section]

\newcommand{\del}{\partial}

\newcommand{\F}{\mathfrak{F}}
\newcommand{\FF}{\mathcal{F}}
\newcommand{\G}{\left|G\right|}
\newcommand{\Z}{{\mathbb{Z}}}
\newcommand{\C}{{\mathbb{C}}}
\newcommand{\R}{{\mathbb{R}}}
\renewcommand{\P}{{\mathfrak{P}}}
\renewcommand{\H}{{\mathbb{H}}}
\newcommand{\Zp}{\mathbb{Z}/p\mathbb{Z}}

\newcommand{\Zz}{\mathbb{Z}/2\mathbb{Z}}

\renewcommand{\i}{\bm{i}}
\renewcommand{\j}{\bm{j}}

\newcommand{\co}{\colon\thinspace}

\begin{document} 

\title{Equivariant Framings, lens spaces and Contact Structures}

\begin{abstract}
We construct a simple topological invariant of certain $3$-manifolds,
including quotients of $S^3$ by finite groups, based on
the fact that the tangent bundle of an orientable $3$-manifold is
trivialisable. This invariant is strong enough to yield the
classification of lens spaces of odd, prime order. We also use
properties of this invariant to show that there is an oriented
$3$-manifold with no universally tight contact structure. We
generalise and sharpen this invariant to an invariant of a finite
covering of a $3$-manifold.
\end{abstract}

\subjclass{Primary: 57M60; 57M27}

\date{\today}
\author{Siddhartha Gadgil}

\address{	Department of mathematics\\
		SUNY at Stony Brook\\
		Stony Brook, NY 11794}
\email{gadgil@math.sunysb.edu}

\maketitle 

It is well known that the tangent bundle of an orientable $3$-manifold
is trivialisable. This is in particular true for manifolds of the form
$M=S^3/G$, where $G$ is a finite group acting without fixed points on
$S^3$. These are the so-called topological spherical space forms.

Using the fact that $TM$ is trivialisable, we can define an invariant
$\F(M)$ of $M$ with a fixed orientation, which we call the
\emph{equivariant framing} of $M$. Namely, the homotopy classes of
trivialisations of the tangent bundle of $S^3$ are a torseur of $\Z$
(i.e., a set on which $\Z$ acts freely and transitively), which can
moreover be canonically identified with $\Z$ by using the Lie group
structure of $S^3$ as the unit quaternions and identifying a
left-invariant framing with $0\in\Z$. Now, find a trivialisation
$\tau$ of $TM$ and pull it back to one of $TS^3$. Under the above
identification, this gives an element $\F(M,\tau)\in\Z$. This
certainly depends on $\tau$, but we shall see that its reduction
modulo $\G$, when $H_1(M,\Z_2)=0$ (in particular when $\G$ is odd),
and modulo $\G/2$ otherwise, is well-defined. Observe that this is the
same as the collection of homotopy classes of \emph{equivariant}
trivialisations with respect to the action of $G$ on $S^3$.

The above definition does not depend on the identification of $S^3$
with the universal cover of $M$, since two such identifications differ
by an orientation preserving self-homeomorphism of $S^3$, which must be
isotopic to the identity.

Notice that this definition makes essential use of the fact that we
have a quotient of $S^3$, rather than a homology (or even homotopy)
sphere. In the more general situation, where we have a quotient of a
homology sphere by a finite cyclic group, we can use canonical
$2$-framings, as introduced by Atiyah in~\cite{At}. This in fact gives
an invariant corresponding to any finite cover of any manifold.

More interestingly, we can obtain an integer-valued invariant. To do
this we compare the pullback of the canonical $2$-framing of $M$ with
that of a finite cover of $M$. We shall define this in
section~\ref{S:general}.

The invariant $\F(L(p,q))$ can readily be computed for odd $p$. As it
turns out, it suffices to classify lens spaces with $p$ prime.

Thus, $\F(\cdot)$ is an invariant sensitive enough to distinguish
between homotopy equivalent manifolds. It is arguably one of the
simplest such invariants.

Furthermore, there is a transparent relation between $\F(\cdot)$ and
the exceptional isomorphism $SO(4)=(SU(2)\times SU(2))/{\pm 1}$. This
means $\F(\cdot)$ is likely to be useful in studying free finite group
actions on $S^3$ - the most elegant classification of orthogonal
actions can be obtained using the exceptional isomorphism, so it is
useful to have a related topological invariant. At the time of writing
this work is in progress.

Another immediate consequence of the existence of this invariant is
that corresponding to one of the orientations of the Poincar\'e
homology sphere, we do not have a positive universally tight contact
structure. Further applications in a similar vein are given
in~\cite{Ga}.

For this application, one can use tangent plane fields that are
trivialisable rather than framings. The framing invariant in this
context is equivalent to an invariant of tangent plane fields defined
by Gompf. The relation between these is explained in
section~\ref{S:Gompf}. Under this equivalence, our proof translates to
Gompf's proof.

Both framings of and tangent plane fields in $3$-manifolds have been
studied classically. There has also been recent work on these,
motivated by relations to Topological Quantum field theories,
Seiberg-Witten invariants and contact geometry. Motivated by the work
of Witten~\cite{Wi}, framings have been studied by Atiyah~\cite{At},
Freed and Gompf~\cite{FG}, Reshetikhin and Turaev~\cite{RT} and Kirby
and Melvin~\cite{KM}. Tangent plane fields on a $3$-manifold were
first classified in terms of a framing by Pontrjagin~\cite{Po}. More
recently, intrinsic invariants of these have been studied and used by
Lisca and Matic~\cite{LM}, G.~Kuperberg~\cite{Go} and finally by
Gompf~\cite{Go}. We clarify the relation with Gompf's work in
section~\ref{S:Gompf}.

The principle novelty of this paper is the use of framings to define a
useful invariant of $3$-manifolds. While Gompf studies the behaviour
under pullbacks of his invariants of tangent plane fields, which is
equivalent for some applications, he does not define or use an
invariant of $3$-manifolds.

In section ~\ref{S:defn} we show that the invariant $\F(M)$ can indeed
be defined as above. In section~\ref{S:lens} it is computed for
odd-order lens spaces. In section~\ref{S:contact} we relate this to
contact structures and prove the result regarding contact
structures. We generalise and sharpen the framing invariant in
section~\ref{S:general}

\begin{acknowledgements} 
I would like to Slava Matveyev and Kalyan Mukherjea for several
helpful conversations and comments. I would also like to thank the
referee for several helpful comments and corrections.
\end{acknowledgements}

\section{The Definition of the Framing invariant}\label{S:defn}

To define the framing invariant, we need some (straightforward)
results.

\begin{proposition} The set of homotopy classes of trivialisations of
$TS^3$ is a torseur of $\Z$.  
\end{proposition} 
\begin{proof} Given two trivialisations, expressing one in terms of the
other gives a map from $S^3$ to $SO(3)$. The homotopy classes of
trivialisations are the homotopy classes of such maps. But as $S^3$ is
simply connected such maps lift to maps $S^3\to SU(2)\cong S^3$, as
do homotopies between them. The homotopy class of a map from $S^3$ to
itself is determined by its degree (as orientations have been fixed),
making the homotopy classes of trivialisations a torseur of $\Z$.
\end{proof}

\begin{definition} Consider $S^3\cong SU(2)$ as the Lie group of unit
quaternions. The \emph{canonical framing} of $TS^3$ is the framing of
$S^3$ which is invariant under left multiplication and is $(i,j,k)$ at
the identity.
\end{definition}

\begin{proposition} The homotopy class of the canonical framing is
determined by an orientation of $S^3$, and does not depend on the
identification with $SU(2)$.
\end{proposition}

\begin{proof} Suppose $f\co S^3\to SU(2)$ is an isomorphism giving a
second Lie group structure to $S^3$. Then we have an induced
orientation-preserving diffeomorphism $\phi\co S^3\to S^3$, and we need
to show that the pullback of the canonical trivialisation under $\phi$
is homotopic to the canonical trivialisation. But it is well known
that any orientation preserving homeomorphism from $S^3$ to itself is
isotopic to the identity. Thus, the pullback fixes the homotopy class
of any trivialisation, in particular the identity trivialisation.
\end{proof}

For the remainder of the section, let $M=S^3/G$, where $G$ is a finite
group that acts without fixed points on $S^3$. Let $\pi\co S^3\to M$
be the projection map.

\begin{proposition} The set of homotopy classes of trivialisations of
$TM$ is a torseur of $\Z$ when $H_1(M,\Zz)=0$(in particular when $\G$
is odd). When $\H_1(M,\Zz)\neq 0$, the trivialisation is determined by
a map $M\to SO(3)$
\end{proposition}

\begin{proof} As with $S^3$, the difference between trivialisations is
determined by a map $f\co M\to SO(3)$. When $H_1(M,\Zz)=0$, this lifts to
a map $\phi\co \pi_1(M)\to S^3$. It is well known that the homotopy class
of such a map is determined by its degree.
\end{proof}

Thus, when $H_1(M,\Zz)=0$, the difference between two trivialisations
of $M$ can be regarded as an integer. In particular this is true for
$S^3$.

\begin{proposition} 
Suppose $\tau_i,i=1,2$ are trivialisations of $TM$ and $\pi^*(\tau_i)$
are their pullbacks.  Then $\pi^*(\tau_1)-\pi^*(\tau_2)$ is divisible
by $\G$ when $H_1(M,\Zz)=0$ and by $\G/2$ when $H_1(M,\Zz)\neq 0$.
\end{proposition}

\begin{proof} 
As above, when $H_1(M,\Zz)=0$, we have a map $\phi\co M\to S^3$
representing the difference between the $\tau_i$. It is easy to see
that the map $\phi\circ\pi$ represents the difference between the
pullbacks. As $\pi$ has degree $\G$ and the degree is multiplicative,
it follows that $\pi^*(\tau_1)-\pi^*(\tau_2)=deg(\phi\circ\pi)$ is
divisible by $\G$.

In the case when $H_1(M,\Zz)\neq 0$ (as also when $H_1(M,\Zz)=0$),
there is a map $\phi\co M\to SO(3)$ representing the difference between
the trivialisations. On composing with the covering map, this gives
the a map representing $\pi^*(\tau_1)-\pi^*(\tau_2)$ which lifts to
$\tilde\phi\co S^3\to S^3$. Thus, if $\alpha:S^3\to SO(3)$ is the
covering map, we get a commutative diagram

$$\begin{CD}
S^3             @>\tilde\phi>>  S^3       \\
@V{\pi}VV	           @VV{\alpha}V\\
M               @>\phi>>       SO(3)
\end{CD}$$

As the degree of maps is multiplicative, and $deg(\pi)=\G$ and
$deg(\alpha)=2$, we get $2\cdot deg(\tilde\phi)=\G\cdot deg(\phi)$, or
$deg(\tilde\phi)=\frac{\G}{2}deg(\phi)$. The result follows

\end{proof} 

We are now in a position to define the invariant $\mathfrak{F}(M)$

\begin{definition} Let $M=S^3/G$, where $G$ is an finite
group acting without fixed points on $S^3$. The framing invariant
$\F(M)\in \Z/\left<G\right>\Z$, where $\left<G\right>=\G$ when
$H_1(M,\Zz)=0$ and $\left<G\right>=\G/2$ when $H_1(M,\Zz)\neq 0$, is
the equivalence class of the trivialisation of $TS^3$ obtained by
pulling back a trivialisation of $M$.
\end{definition}

\section{Computation for lens spaces}\label{S:lens} 

It is easy to see that $\F(\cdot)$ is a non-trivial invariant, and is
in fact sensitive enough to distinguish between homotopy equivalent
manifolds. For, the lens spaces $L(p,1)$ are quotients of $S^3$ by a
subgroup of the unit quaternions acting on themselves by left
multiplication, and hence the canonical trivialisation is
equivariant. On the other hand, the spaces $L(p,-1)$ are quotients by
a subgroup of the unit quaternions acting on themselves by right
multiplication, making the right-invariant trivialisation
equivariant. These two trivialisations differ by the adjoint action of
$SU(2)$ on its Lie algebra. This lifts to a degree $1$ map from
$SU(2)$ to itself. Thus the framing invariant suffices to show that
there is no orientation preserving homeomorphism between $L(p,1)$ and
$L(p,-1)$ when $p\neq 2$.

But for $p\equiv 1 (mod\ 4)$, with $p$ a prime, there is an
orientation preserving homotopy equivalence between $L(p,1)$ and
$L(p,-1)$. By uniqueness of prime decompositions of $3$-manifolds, it
follows, for example, that $L(5,1)\#L(5,1)$ is homotopy equivalent but
not homeomorphic to $L(5,1)\#L(5,-1)$.

Thus, $\F(\cdot)$ depends essentially on the homeomorphism type, and
not just the homotopy type, of a lens space. Our goal here is to
compute this explicitly for odd-order lens spaces, and show that if
$p$ is an odd prime, $\F(\cdot)$ suffices to classify lens spaces.

\begin{theorem} Suppose $p$ is odd. Then
$\F(L(p,q))=\frac{(q-1)(q^{-1}-1)}{4}$, where $q^{-1}$ is a
multiplicative inverse of $q$ modulo $p$, and $q$ and $q^{-1}$ have
been chosen to be odd representatives of their mod $p$ equivalence
class.  
\end{theorem} 

\begin{remark} The formula above does not depend on the choice of the odd
representatives for $q$ and $q^{-1}$
\end{remark}

\begin{proof}

It will be useful to regard $S^3$ as the join $S^1*S^1$, which is
embedded in the natural way in the quaternions. Then the action
corresponding to the lens space $L(p,q)$ is the join of actions on the
circle $S^1\subset \C$ generated respectively by $z_1\mapsto
e^{\frac{2\pi \i}{p}}z_1$ and $z_2\mapsto e^{\frac{2\pi \i q}{p}}z_2$. We
first find an equivariant trivialisation along the circles $(z_1,0)$
and $(0,z_2)$ and then extend these to $S^3$. 

Henceforth the equivariant trivialisation is given in terms of the map
$f\co S^3\to SU(2)\cong S^3$ representing the difference with the
left-invariant trivialisation. Along the circle $C_1=\{(z_1,0)\co z_1\in
S^1\}$, we can take the first vector to point along the circle. An
equivariant trivialisation is obtained by choosing the other two
vectors so that they rotate $q$ times, for some choice in the mod $p$
class of $q$. This follows as the arc joining $(1,0)$ to
$(e^{\frac{2\pi \i}{p}},0)$ is a fundamental domain, with the only
identification due to the group action being that of the endpoints
induced by the element $(z_1,z_2)\mapsto (e^{\frac{2\pi
\i}{p}}z_1,e^{\frac{2\pi \i q}{p}}z_2)$, and this element rotates the
normal plane to the circle by $e^{\frac{2\pi \i q}{p}}$.

This trivialisation differs from the left-invariant trivialisation by
$q-1$ rotations, and if $q$ is chosen odd, this gives a comparison map
which lifts to $S^3$. The resulting restriction of $f$ is the map
$f\co (z_1,0)\mapsto (z_1^{\frac{q-1}{2}},0)$ of degree $\frac{q-1}{2}$
from the circle to itself.

Likewise, we can trivialise the tangent bundle along the circle
$C_2=\{(0,z_2)\co z_2\in S^1\}$. It is convenient to choose the
trivialisation so that it differs from the left invariant one at
$(0,1)$ by $\j=(0,1)$. Then the map $f$ restricted to this circle is
the map $f\co (0,z_1)\mapsto (0,z_2^{\frac{q^{-1}-1}{2}})$ of degree
$\frac{q^{-1}-1}{2}$ from the circle to itself as the identification
on the boundary of the fundamental domain is induced in this case by
$(z_1,z_2)\mapsto (e^{\frac{2\pi iq^{-1}}{p}}z_1,e^{\frac{2\pi
\i}{p}}z_2)$.

We now extend this map to a map from the disc $D_0=\{(z_1,r)\in
S^3\co \left|z_1\right|\leq 1,r\in\R\}\subset S^3\subset\C$ of the form
$re^{i\theta}\mapsto re^{i\frac{q-1}{2}\theta},r,\theta\in \R$. Since
the only points in the disc which are identified by the group action
are on the boundary, and the trivialisation on the boundary is
equivariant, we get an equivariant trivialisation of the disc
corresponding to this map.

We extend this by requiring equivariance to the disc
$D_1=\{(z_1,re^{\frac{2\pi \i}{p}})\in S^3:  \left|z_1\right|\leq
1,r\in \R\}$. The map is previously defined on the boundary of this
disc, where it is equivariant. Thus, we have a uniquely equivariant
extension.

The midpoint of $D_1$ is part of the circle $C_2$, as is an arc
$\alpha$ joining the midpoints of $D_1$ and $D_2$. The map $f$ as
previously defined on $C_2$, must agree with the definition on $D_1$
at the midpoint as both these have been defined so that the
trivialisation is equivariant.

The discs $D_0$ and $D_1$ bound a \emph{lens} $N$, which is a
fundamental domain for the group action. The map $f$ has already been
defined on the boundary discs as well as the arc $\alpha$. This
extends to a map on $N$, and hence a trivialisation of the tangent
bundle. Up to homotopy, any other choice of $f$ can be obtained by
replacing $f$ in a small neighbourhood of an interior point $x$ of the
fundamental domain by a degree $k$ map from $S^3$ to itself for some
$k$. More precisely, $f$ is homotopic to a map that is constant on the
neighbourhood of the point $x$. Replace this map in this neighbourhood
by a map to $S^3$ that maps the boundary of the neighbourhood to a
single point and so that the inverse image of a generic point has
algebraic multiplicity $k$. We shall call such a transformation
`blowing a degree--$k$ bubble'.

We now extend the trivialisation equivariantly from the fundamental
domain to all of the manifold, and define $f$ accordingly. This is an
extension of the map on $C_1$ and $C_2$ that was previously defined.

The resulting map $f\co S^3\to S^3$ is homotopic to a map obtained from
the join of maps $C_i\to C_i$ defined by taking powers on the unit
circle in $\C$ by blowing a degree--$k$ bubble in each of the $p$
images of the fundamental domain. Thus, it has degree
$\frac{(q-1)(q^{-1}-1)}{4}+kp$. This proves our claim.
\end{proof}

\begin{corollary} 
Suppose $p$ is a prime, then $L(p,q)=L(p,q')$ as oriented manifolds 
if and only if $q'=q^{\pm 1}$
\end{corollary}
\begin{proof}
It is well known that $L(p,q)=L(p,q^{\pm 1})$ as oriented manifolds (a
homeomorphism is induced by $(z_1,z_2)\mapsto (z_2,z_1)$. Conversely,
for a fixed $q$, the condition $\F(L(p,q'))=\F(L(p,q))$ is a quadratic
equation in $q'$ over the field $\Zp$, with roots $q$ and $q^{-1}$.

If $q^{-1}\neq q$, then these are two distinct root of the quadratic
equation, and hence the only solutions for $q'$. If $q=q^{-1}$, and
$q'\neq q$ is another root, then we also have $q'=(q')^{-1}$,
otherwise we would have three distinct roots. But this means that
$q'=\pm 1$ and $q=\pm 1$, and we have already seen that
$\F(L(p,1))\neq\F(L(p,-1))$. Alternatively, using $q'=q^{-1}$, the
above reduces to a linear equation for $q'$ that is satisfied by $q$,
and hence $q'=q$.
\end{proof}

The same statement is well known to be true for all values of $p$ and
there are several proofs of this (see, for instance, \cite{Re},
\cite{Br}, \cite{BO} and \cite{RS}).

\begin{remark}
It is more natural to declare the left invariant framing to be
$-\frac{1}{2}$ rather than $0$. Then we have the relation
$\F(L(p,-q))=-\F(L(p,q))$.
\end{remark}

\begin{remark}
More generally, after re-normalising as above, $\F(-M)=-\F(M)$. This
is immediate from section~\ref{S:general} and can also be proved
directly.
\end{remark}

\section{Universally tight contact structures}\label{S:contact}

Let $M$ be a closed, orientable $3$-manifold. Recall that a contact
structure $\xi$ on $M$ is a totally non-integrable tangent plane
field. We shall assume that the tangent plane field is co-orientable
(we shall say that the contact structure is co-orientable). In this
situation, we can express $\xi=ker(\alpha)$, where $\alpha$ is a
$1$-form consistent with the co-orientation.

The hypothesis of $\xi$ being nowhere integrable is equivalent to
$\alpha\wedge d\alpha$ being a non-degenerate $3$-form. Thus, this is
everywhere a non-zero multiple of the volume form, and hence induces
an orientation on $M$. We say that $\xi$ is \emph{positive} if this
orientation agrees with the orientation of $M$.

A fundamental dichotomy among contact structures on $3$-manifolds is
between \emph{tight} and \emph{overtwisted} contact structures. An
\emph{overtwisted} contact structure is one that contains an unknot
that is everywhere tangent to the contact structure so that the
framing induced by the contact structure is the $0$-framing. A contact
structure that is not overtwisted is said to be \emph{tight}. A
\emph{universally tight} contact structure is one that pulls back to
a tight contact structure on every cover of $M$.

A fundamental result of Eliashberg is that $S^3$ with a fixed
orientation has a unique positive tight contact structure, namely the
contact structure invariant under left multiplication. Our results
follow from this and some simple observations.

\begin{proposition}~\label{T:Fr} 
Let $M$ be an integral homology $3$-sphere with a contact structure
$\xi$. Then there is a canonical framing associated to $\xi$. Further,
the pullback of this framing to any homology sphere that covers $M$ is
the framing induced by the pullback of the contact structure.
\end{proposition}
\begin{proof} 
As $M$ is a homology sphere, the contact-structure is
co-orientable. Choose and fix a co-orientation. This induces an
orientation on the plane-bundle given by the contact structure, which
we identify with $\xi$.

As $H^2(M)=0$, the Euler class of $\xi$ is trivial. Hence there is a
trivialisation of $\xi$. Further, two trivialisations differ by a map
onto $S^1$. As $H^1(M)=0$, any such map is homotopic to a constant
map.

Thus, there is a trivialisation $(X_1,X_2)$ of $\xi$, canonical up to
homotopy. This, together with a vector $X_3$ normal to $\xi$, that is
consistent with the co-orientation, gives a framing $(X_1,X_2,X_3)$.

The homotopy class of this trivialisation does not depend on the
choice of co-orientation since $(X_1,-X_2,-X_3)$ gives a
trivialisation corresponding to the opposite co-orientation, and this
is clearly homotopic to the trivialisation $(X_1,X_2,X_3)$.

As the trivialisation of $\xi$ pulls back to give a trivialisation of
the pullback to any cover of $\xi$, the second claim follows.
\end{proof}

Now let $\P$ be the Poincare homology sphere with a fixed orientation,
and let $-\P$ denote the same manifold with the opposite
orientation. These manifolds have finite fundamental group. The
Poincar\'e homology sphere is the quotient of $S^3$ by a group acting
by left multiplication, and $-\P$ is the quotient of an action by
right multiplication. We can now prove the following theorem. Note
that any contact structure on a homology sphere is automatically
co-orientable.

\begin{theorem}[Gompf]~\label{T:tght}
The manifold $-\P$ does not have a universally tight positive contact
structure.
\end{theorem} 
\begin{proof}
As $-\P$ is the quotient of $S^3$ by a group acting by right
multiplication, it follows that any framing on $-\P$ pulls back to one
homotopic to a framing invariant under right Lie multiplication, or
one differing from this by $\left|\pi_1(\P)\right|$ units (as
$H_1(\P,\Z_2)=0$). However, if $-\P$ had a universally tight positive
contact structure, then the associated framing must pulls back to give
the framing associated to left Lie multiplication. This gives the
required contradiction.
\end{proof}

\begin{remark} 
Etnyre and Honda~\cite{EH} have shown that $-\P$ does not have a tight
contact structure.
\end{remark}

\begin{remark} 
We see in section~\ref{S:Gompf} that the proof of the above result
translates to Gompf's proof under the correspondence between framings
and trivialisable tangent plane fields.
\end{remark}

V.Colin~\cite{Co} shows that tight contact structures on connected
sums of manifolds are connected sums of tight contact structures on
each summand. Using this, we obtain the following result.

\begin{corollary} The manifold $\P \#-\P$ does not admit a universally tight
contact structure.
\end{corollary}


\section{2-Framings and invariants of covers}\label{S:general}

We now generalise and sharpen the framing invariant using so called
canonical $2$-framings as introduced by Atiyah~\cite{At}. Atiyah has
shown that any $3$-manifold has associated to it a canonical framing
$\FF$ of the Whitney sum $2TM=TM\oplus TM$ of the tangent bundle with
itself, considered as a Spin(6)-bundle with the natural spin
structure. The framing $\FF$ is characterised by

$$\sigma(W^4)=\frac{1}{6}p_1(2TW,\mathcal{F})$$
for any smooth $4$-manifold $W$ with $\del W=M$. Here $\sigma$
denotes the signature and $p_1$ the relative Pontrjagin class. By the
Hirzebruch signature formula this does not depend on the choice of $W$.

Atiyah has shown that such a $2$-framing always exists, and the
$2$-framings form a torseur of $\Z$. Now, we can define the framing
invariant $\F(M,N)$ associated to a cover $M\to N$ - pull back the
canonical $2$-framing of $N$ and compare this with the canonical
$2$-framing of $M$. Thus, we get an integer-valued invariant. We state
for reference the following lemma, which is an immediate consequence
of Atiyah's result.

\begin{lemma}\label{T:defect}
Suppose $M$ is a $3$-manifold bounding a $4$-manifold $W$ and let
$\FF$ be any $2$-framing of $M$. The difference between $\FF$ and the
canonical framing is $\sigma(W^4)-\frac{1}{6}p_1(2TW,\mathcal{F})$
\end{lemma}

As a framing gives a $2$-framing, we see that we have a sharpening of
the framing invariant defined earlier. We show here that this is a
non-trivial invariant. 

\begin{remark}
A hyperbolic $3$-manifold has many covers, hence many invariants
associated to it. It is not clear whether these are useful.
\end{remark}

\begin{theorem} 
Suppose $N$ and $N'$ are $h$-cobordant $3$-manifolds, $M$ is a cover
of $N$ and $M'$ the corresponding cover of $N'$. Then $\F(M,N)=\F(M',N')$.
\end{theorem}
\begin{proof}
Let $X$ be the $h$-cobordism between $N$ and $N'$, so that $\del
W=N'-N$. If $W$ is a $4$-manifold with boundary $N$, then
$W'=W\coprod_{N}X$ is a $4$-manifold with boundary $N'$ with the same
signature as $W$. By considering these manifolds, it is immediate that
that if $\F$ is the canonical $2$-framing for $M$, then the canonical
framing $\F'$ of $N'$ is characterised by 
$$p_1(2TX,\F,\F')=0$$

The $h$-cobordism $X$ lifts to an $h$-cobordism $Y$ between $M$ and
$M'$, and the Pontrjagin class relative to the framings pulled back is
the pullback of the Pontrjagin class and hence is zero. Let $U$ be a
$4$-manifold with boundary $M$, and let $U'=U\coprod_{M}Y$. Applying
lemma~\ref{T:defect} to $U$ and to $U'$, the result follows.
\end{proof}

We can generalise the above theorem to the following.

\begin{theorem} 
Let $X$ be a cobordism between $N_1$ and $N_2$ and let
$\phi\co \pi_1(W)\to H$ be a surjective map onto a finite group that
restricts to surjections on $\pi_1(N_1)$ and $\pi_1(N_2)$. Suppose
that the cover $\tilde X$ with fundamental group $ker(\phi)$ satisfies

$$\sigma(\tilde W)=\left| H\right|\sigma(W)$$

and $M_1$ and $M_2$ are the covers of $N_i$ with fundamental group
$ker(\phi)$. Then $\F(M,N)=\F(M',N')$.
\end{theorem}
\begin{proof}
We use the additivity of the signature and relative Pontrjagin
classes. The above proof generalises immediately.
\end{proof} 

\begin{example} 
Let $N_1$ be a lens space, $M_1$ be $S^3$ and $K$ be a homologically
trivial knot in $N_1$ the components of whose lift to $S^3$ are
unlinked (for instance, the untwisted Whitehead double of any
homologically trivial knot). Add a $2$-handle to $N_1$ along $K$ with
framing $1$ to get a $4$-manifold $W$, and let its other boundary
component be $N_2$. Let $\phi$ be the map onto $\pi_1(N_1)$ that
extends the identity map on $N_1$. Then $\F(M,N)=\F(M',N')$.
\end{example}

The manifold $N_2$ is the result of surgery about the knot $K$. Using
this construction and taking connected sums, we get hyperbolic
manifolds with covers having various framings.

\section{Relation to Gompf's invariants}~\label{S:Gompf}

In this section we relate equivariant framings to Gompf's invariant,
equivalent to the $3$-dimensional obstruction in Pontrjagin's
classification, for tangent plane fields. Assume henceforth that $M$
is an oriented rational homology sphere (for instance, $M=S^3/G$). We
first establish a canonical correspondence between homotopy
classes of framings $\FF$ of $M$ and homotopy classes of orientable
tangent plane fields $\xi$ on $M$ with $c_1(\xi)=0$.

\begin{proposition}
Let $M$ be an rational homology $3$-sphere. Then there is a natural
bijective correspondence between homotopy classes of framings $\FF$ of
$M$ and homotopy classes of orientable tangent plane fields $\xi$ on
$M$ with $c_1(\xi)=0$.
\end{proposition}
\begin{proof}
We proceed as in Proposition~\ref{T:Fr}. Given an orientable tangent
plane fields $\xi$ on $M$ with $c_1(\xi)=0$, fix an orientation on
$\xi$. There is a trivialisation $(X_1,X_2)$ of $\xi$ respecting the
orientation. Further, the homotopy classes of such trivialisations are
classified by $H^1(M)=0$. Hence, as $M$ is a rational homology sphere,
the trivialisation is canonical up to homotopy.

Find a vector field $X_3$ normal to $\xi$ such that $(X_1,X_2,X_3)$
respects the orientation on $M$. This is possible as $\xi$ and $M$ are
orientable. Then $\FF=(X_1,X_2,X_3)$ is a framing of $M$. As before
this does not depend on the choice of orientation.

Conversely, given a framing $\FF=(X_1,X_2,X_3)$ let $\xi$ be the span
of $X_1$ and $X_2$. Evidently these constructions are inverses of each
other.
\end{proof}

We now recall the invariant of Gompf for such tangent plane fields. To
do this one finds an almost complex $4$-manifold $(X,J)$ with $\del
X=M$ so that the plane field $TM\cap J(TM)$ on $M$ induced by the
almost complex structure is $\xi$. Using this, one would wish to
define the invariant $c_1^2(X)-2\chi(X)-3\sigma(X)$, where $\chi$ and
$\sigma$ denote the Euler characteristic and the signature.

One cannot always define such an invariant as $c_1(X)\in
H^2(X)\cong H_2(X,\del X)$ and there is no natural pairing on $H_2(X,\del
X)$. Gompf works instead with a pairing on surfaces representing
elements of $H_2(X,\del X)$ with given framings on their boundaries.

In case of tangent plane fields $\xi$ on $M$ with $c_1(\xi)=0$, we can
define the invariant directly by the following lemma.

\begin{lemma}
Let $X$ be an almost complex $4$-manifold $X$ with $\del X=M$ so that
the plane field $\xi$ on $M$ induced by the almost complex structure
satisfies $c_1(\xi)=0$. Then the Poincar\'e dual $PD(c_1(X))\in
H_2(X,\del X)$ of $c_1(X)$ is contained in the image of $H_2(X)$ under
the inclusion map.
\end{lemma}
\begin{proof}
If $[F]=PD(c_1(X))$ for a properly embedded surface $F$ with boundary,
then it is easy to see that $[\del F]=PD(c_1(\xi))$ as $TX|_{\del X}$
splits as the sum of a trivial complex line bundle and $\xi$. It
follows that $[\del F]=0$. Hence $PD(c_1(X))$ is in the image of
$H_2(X)$ by the long exact sequence of homology groups.
\end{proof}

Thus, Gompf's invariant reduces to
$c_1^2(X)-2\chi(X)-3\sigma(X)$. Using $p_1(X,\FF)=c_1^2(X)+2c_2(X)$
and $c_2(X)=\chi(X)$, we see that this is
$-3(\sigma(X)-\frac{1}{6}p_1(2TX,\FF))$. We saw that
$\sigma(X)-\frac{1}{6}p_1(2TX,\F)=0$ characterises the canonical
framing, which is used to give an integral sharpening of the framing
invariant.

To prove theorem~\ref{T:tght}, Gompf uses the fact that if we pullback
two tangent plane fields from $-\P$ to $S^3$, then the difference
between the values of the invariant for the two plane fields is a
multiple of $4$ times the degree of the cover (the factor of $4$ comes
about in the process of defining the invariant for general plane
fields). This is equivalent to our argument using equivariant
framings.


\begin{thebibliography}{10}

\bibitem{At} M.Atiyah 
\textit{On framings of $3$-manifolds} 
Topology \textbf{29} (1990), 1--7. 

\bibitem{BO} F.Bonahon and J.P.Otal
\textit{Sciendements de Heegard des espaces lenticulaires}
Ann. Sci. \'Ecole Norm. Sup. \textbf{16} (1983), 451--466.

\bibitem{Br} E.J.Brody
\textit{The topological classification of lens spaces}
Annals of mathematics \textbf{71} (1960), 163--184.

\bibitem{Co} V.Colin 
\textit{Chirurgies d'indice un et isotopies de sphères
dans les variétés de contact tendues} 
C. R. Acad. Sci. Paris Sér. I Math. \textbf{324} (1997), 659--663

\bibitem{El} Y.Eliashberg 
\textit{Contact $3$-manifolds twenty years since J. Martinet's
work}.
Ann. Inst. Fourier (Grenoble) \textbf{42} (1992), 165--192.


\bibitem{EH} J.Etnyre and K.Honda
\textit{On the non-existence of tight contact structures}
Ann. of Math. (2) \textbf{153} (2001), 749--766. 

\bibitem{Ga} S.Gadgil,
\textit{Contact structures on elliptic manifolds}
preprint, 2001.

\bibitem{FG} D.Freed and R.Gompf 
\textit{Computer calculation of Witten's $3$-manifold invariant.} 
Comm. Math. Phys. \textbf{141} (1991), no. 1, 79--117

\bibitem{Go} R.Gompf
\textit{Handlebody construction of Stein surfaces}
Ann. of Math. (2) \textbf{148} (1998), 619--693.
 

\bibitem{KM} R.Kirby and P.Melvin
\textit{Canonical framings for $3$-manifolds.} 
Turkish J. Math. \textbf{23} (1999), no. 1, 89--115. 


\bibitem{Ku} G.Kuperberg 
\textit{Noninvolutory Hopf algebras and $3$-manifold invariants.} 
Duke Math. J. \textbf{84} (1996), 83--129

\bibitem{LM} P.Lisca and G.Matic
\textit{On homotopic, non-isomorphic tight contact structures 
on $3$-manifolds.} 
Turkish J. Math. \textbf{20} (1996), no. 1, 37--45. 

\bibitem{Po} L.Pontrjagin 
\textit{A classification of mappings of the three-dimensional 
complex into the two-dimensional sphere.} 
Rec. Math. [Mat. Sbornik] N. S. \textbf{9 (51)}, (1941). 331--363

\bibitem{Re} K.Reidemeister
\textit{Homotopieringe und Linsenr\"aume}
Abh. Math. Sem. Hamburg \textbf{11} (1935), 102--109.

\bibitem{RT} Reshetikhin, N. and Turaev, V. G.  
\textit{Invariants of $3$-manifolds via link polynomials 
and quantum groups.}
Invent. Math. \textbf{103} (1991) 547--597.


\bibitem{RS} J.Hyam Rubinstein and M.G.Scharlemann
\textit{Comparing Heegard splittings of non-Haken $3$-manifolds}
Topology \textbf{35} (1996) 1005--1026.

\bibitem{Wi} E.Witten
\textit{Quantum field theory and the Jones polynomial.} 
Comm. Math. Phys. \textbf{121} (1989), no. 3, 351--399. 

\end{thebibliography}
\end{document}